\documentclass{article}[10pt]
\usepackage{amssymb}

\font\twelvebb=msbm10
\font\fivebb=msbm5
\font\sevenbb=msbm7
\newfam\bbfam  \scriptscriptfont\bbfam=\fivebb
\textfont\bbfam=\twelvebb
\scriptfont\bbfam=\sevenbb
\def\bb{\fam\bbfam\twelvebb}

\def\N{{\bb N}}
\def\Z{{\bb Z}}
\def\R{{\bb R}}

\catcode`\é=\active
\defé{\'e}
\catcode`\è=\active
\defè{\`e}
\catcode`\ê=\active
\defê{\^e}
\catcode`\à=\active
\defà{\`a}
\catcode`\â=\active
\defâ{\^a}
\catcode`\ù=\active
\defù{\`u}
\catcode`\ô=\active
\defô{\^o}

\def\tab{\hspace*{6mm}}

\newcommand{\fin}{
\vskip 2mm
\noindent
\hfill{--------}}

\def\tab{\hspace*{6mm}}
\begin{document}
\begin{center}
\LARGE Essentialité dans les bases additives
\end{center}
\vspace{-3cm}
\begin{flushleft}
{\it J. Number Theory}, \\
{\bf 123} (2007), p. 170-192.
\end{flushleft}~\vspace{1cm}
\noindent \begin{center} BRUNO DESCHAMPS
(Bruno.Deschamps@univ-lemans.fr) \\ et \\ BAKIR FARHI
(bakir.farhi@gmail.com) ~\vspace{2mm} \\
Département de Mathématiques, Université du Maine \\
Avenue Olivier Messiaen, 72085 Le Mans Cedex 9 \\ France
\end{center} \vskip 2mm \noindent
{\small {\bf R\'esum\'e.---} Dans cet article nous étudions la notion de partie essentielle d'une base additive,
c'est-à-dire de parties finies minimales $P$ d'une base $A$ donnée telles que $A \setminus P$ ne soit plus une base. L'existence de parties essentielles pour une base équivaut à ce que cette base soit incluse, à partir d'un certain rang, dans une progression arithmétique non triviale. Nous montrons que pour toute base $A$ il existe une progression arithmétique de plus grande raison qui contient, à partir d'un certain rang, la base $A$. Possédant cette plus grande raison $a$, on peut alors majorer le nombre de parties essentielles de $A$ (qui est donc toujours fini) par la longueur du radical de $a$ (i.e. le nombre de nombre premiers divisant $a$).\\
Dans le cas des parties essentielles de cardinal $1$ (éléments
essentiels) nous introduisons une méthode pour "dessentialiser"
une base. En application de ces considérations nous raffinons et
complétons de manière définitive le résultat de Deschamps et
Grekos sur la majoration du nombre d'éléments essentiels dans une
base en fonction de l'ordre : nous montrons que pour toute base
$A$ d'ordre $h$, le nombre $s$ d'éléments essentiels de $A$
vérifie $s\leq c\sqrt{\frac{h}{\log h}}$ avec
$c=30\sqrt{\frac{\log 1564}{1564}}\simeq 2,05728$ et que cette
inégalité est optimale. \vskip 3mm \noindent {\bf Abstract.---} In
this article we study the notion of essential subset of an
additive basis, that is to say the minimal finite subsets $P$ of a
basis $A$ such that $A \setminus P$ doesn't remains a basis. The
existence of an essential subset for a basis is equivalent for
this basis to be included, for almost all elements, in an
arithmetic non-trivial progression. We show that for every basis
$A$ there exists an arithmetic progression with a biggest common
difference containing $A$. Having this common difference $a$ we
are able to give an upper bound to the number of essential subsets
of $A$:
this is the radical's length of $a$ (in particular there is always many finite essential subsets in a basis).\\
In the case of essential subsets of cardinality $1$ (essential
elements) we introduce a way to "dessentialize" a basis. As an
application, we definitively improve the earlier result of
Deschamps and Grekos giving an upper bound of the number of the
essential elements of a basis. More precisely, we show that for
all basis $A$ of order $h$, the number $s$ of essential elements
of $A$ satisfy $s\leq c\sqrt{\frac{h}{\log h}}$ where
$c=30\sqrt{\frac{\log 1564}{1564}}\simeq 2,05728$, and we show
that this inequality is best possible.} \vskip 8mm \noindent
\newpage
\centerline{\bf 1.--- Introduction et notations.} \vskip 2mm
\noindent \tab Dans ce texte, on notera \vskip 2mm \noindent
$\bullet$ $\# E$ le cardinal d'un ensemble $E$, \vskip 2mm
\noindent $\bullet$ $A \setminus B$ (pour deux ensembles $A,B$)
l'ensemble des éléments de $A$ qui ne sont pas dans $B$, \vskip
2mm \noindent $\bullet$ $aX+b$ ($a,b\in \N$ et $X\subset \N$)
l'ensemble d'entiers $\{ax+b ~|~ x\in X\}$, \vskip 2mm \noindent
$\bullet$ $hA$ ($h\in \N^*$ et $A\subset \N$) l'ensemble d'entiers
$\{a_1+\cdots +a_h ~|~ a_1,\cdots ,a_h\in A\}$, \vskip 2mm
\noindent $\bullet$ $E\sim F$ (pour $E,F$ deux parties de $\N$)
pour dire que la différence symétrique de $E$ et de $F$ est finie,
ce qui revient à dire qu'il existe un entier $n_0$ tel que $E\cap
[n_0,+\infty[=F\cap [n_0,+\infty[$, \vskip 2mm \noindent $\bullet$
${\rm ord}(A)$ ($A$ base additive) l'ordre de la base $A$, \vskip
2mm \noindent $\bullet$ $(p_n)_n$ la suite des nombres premiers
énumérés par ordre croissant, \vskip 2mm \noindent $\bullet$
$a_1\cdots \widehat{a_i}\cdots a_n$ ($a_1,\cdots ,a_n\in \N^*$)
l'entier $\displaystyle \frac{a_1\cdots a_n}{a_i}$. \vskip 2mm
\noindent $\bullet$ $\lfloor x \rfloor$ ($x\in \R$) la partie
entière par défaut du réel $x$. \vskip 5mm \noindent \tab Une base
additive (ou base ou encore ensemble de base) est une partie $A$
de $\N$ telle qu'il existe un entier $h\geq 1$ vérifiant $hA\sim
\N$. Le plus petit entier $h$ pour cette propriété s'appelle alors
l'ordre de la base. \vskip 2mm \noindent \tab Dans une base
additive tous les éléments n'ont pas la même importance en terme
de basicité. Par exemple, dans la base $A=2\N\cup \{1\}$ on voit
que $1$ joue un r\^ole prépondérant : si l'on retire $1$ à $A$ on
ne dispose plus d'une base (c'est d'ailleurs dans cet exemple le
seul élément qui ait cette propriété). De manière plus générale
dans une base $A$ il peut y avoir des parties qui sont
essentielles à la basicité, c'est-à-dire des parties $P\subset A$
telles que $A \setminus P$ ne soit plus une base. Si $P$ est une
telle partie et que $P\subset P^{'}$, on voit que $P^{'}$ est
aussi essentielle. Ceci nous amène à la définition suivante :
\vskip 2mm \noindent {\bf Définition.---} {\it Soit $A$ une base
additive. On appelle {\rm essentialité} de $A$ toute partie
$P\subset A$ telle que $A \setminus P$ ne soit pas une base et
telle que $P$ soit minimal au sens de l'inclusion pour cette
propriété. \vskip 2mm \noindent \tab Une essentialité de cardinal
fini sera appelée {\rm partie essentielle} de $A$. Un élément
$a\in A$ tel que $P=\{a\}$ soit une partie essentielle sera appelé
{\rm élément essentiel} de $A$.} \vskip 2mm \noindent \tab Si l'on
considère la base $A=\N$ et si l'on prend un nombre premier $p$ on
voit que la partie $A_p=\N \setminus p\N$ est une essentialité de
$A$. On voit donc qu'il existe des bases additives possédant une
infinité d'essentialités. Dans le cas des éléments essentiels, une
telle chose ne peut pas se produire : Grekos a montré que
l'ensemble des éléments essentiels d'une base était un ensemble
fini (voir [Gr]). Dans [DG] Deschamps et Grekos montrent que ce
nombre est toujours majoré par $5,7\sqrt{\frac{h}{\log h}}$ o\`u
$h$ désigne l'ordre de la base et que cette majoration en
$O\left(\sqrt{\frac{h}{\log h}} \right)$ est la meilleure
possible. Ce résultat est donc optimum, à la constante près. Dans
cet article nous donnons la meilleure constante possible pour
cette estimation : nous montrons (théorème 7) que pour toute base
additive $A$ d'ordre $h$, le nombre $s$ d'éléments essentiels de
$A$ vérifie $s\leq c\sqrt{\frac{h}{\log h}}$ avec
$c=30\sqrt{\frac{\log 1564}{1564}}\simeq 2,05728$ et qu'il existe
une base additive pour laquelle cette majoration est une égalité.
Ceci complète de manière définitive l'étude globale de la
majoration de $s$ en fonction de $h$. \vskip 2mm \noindent \tab
L'étude des parties essentielles d'une base est un prolongement
naturelle de celle des éléments essentiels. Nous montrons dans cet
article (Théorème 11) que, comme pour les éléments essentiels, le
nombre de parties essentielles d'une base est fini. Toutefois, et
au contraire du cas des éléments essentiels, le nombre de parties
essentielles n'est pas majorable par une fonction de l'ordre de la
base. Par exemple, la base additive $X_n=p_1\cdots p_n\N
\cup\{1,2,\cdots ,p_1\cdots p_n\}$ est une base d'ordre $2$ qui
possède $n$ parties essentielles (qui sont les parties
$P_k=\{1\leq i\leq p_1\cdots p_n ~|~ p_k\hskip -1mm\not \hskip
-.13mm |\hskip 1mm i\}$).
Il faut donc, pour controler le nombre de parties essentielles, trouver d'autres invariants susceptibles de donner cette information.\\
\tab En fait l'existence de parties essentielles dans une base additive est à relier à la propriété pour cette base d'être incluse, à partir d'un certain rang, dans une progression arithmétique non triviale (Proposition 8). Dans cet article nous prouvons que pour toute base additive $A$ il existe un plus grand entier $a$ tel que $A$ soit incluse, à partir d'un certain rang, dans une progression arithmétique de raison $a$ (Théorème 9). Ce résultat permet d'introduire trois invariants fondamentaux pour une base additive $A$ : le plus grand entier $a$ donné par le théorème, appelé {\it raison} de la base $A$, et deux parties $B,E\subset \N$ appelées respectivement {\it la dessentialisée} et {\it le réservoir} de $A$ :
les ensembles $B$ et $E$ sont les uniques ensembles tels que $A=(aB+b)\cup E$ ($b\in \{0,\cdots ,a-1\}$)
et tel que pour tout $x\in E$ on ait $x\not\equiv b\ \hbox{\rm mod}(a)$.\\
\tab L'ensemble $B$ est alors une base sans partie essentielle qui
est unique (à translation près et modulo la relation $\sim$) pour
cette propriété (Proposition 10). Ces invariants sont les plus
pertinents pour l'étude des parties essentielles. En effet, nous
montrons qu'étant donnée une base $A$ de raison $a$ et de
réservoir $E$, le nombre de parties essentielles de $A$ est majoré
par la longueur du radical de $a$ (i.e. le nombre de nombres
premiers divisant $a$)  et que toute partie essentielle de $A$ est
incluse dans $E$ (Théorème 11). En fin de texte, un exemple vient
montrer que la majoration du nombre de parties essentielles d'une
base par la longueur du radical de sa raison est la meilleure
possible. \vskip 2mm \noindent \tab Ce texte débute par un
paragraphe qui présente un moyen de dessentialiser élémentairement
une base, c'est-à-dire d'associer de manière naturelle à une base
une autre base sans élément essentiel. Il s'agit en quelque sorte
de l'analogue de la dessentialisée pour le cas des éléments
essentiels. Outre le fait de présenter un algorithme de
dessentialisation, ce paragraphe permet d'introduire des outils et
des résultats qui seront utiles pour les autres parties du texte.
Ils permettent aussi de raffiner des estimations sur l'ordre de
certaines bases (Théorème 3), en particulier ils permettent
d'obtenir la constante optimale pour le théorème de
Deschamps-Grekos. \vskip 10mm \noindent \centerline{\bf 2.---
Eléments essentiels dans une base additive} \vskip 2mm \noindent
\tab Le point de départ pour l'étude des éléments essentiels, et
plus tard des parties essentielles, est une idée de Erdos et
Graham (voir [EG]) exploitée dans l'article [DG] sous la forme de
ce lemme : \vskip 2mm \noindent {\bf Lemme 1.---} {\it Soient $A$
une base additive et $a\in A$. Les propositions suivantes sont
\'equivalentes: \vskip 2mm \noindent i) $A \setminus \{a\}$ est
une base additive, \vskip 2mm \noindent ii) $\hbox{\rm pgcd}
\{b-b^{'} ~|~ b,b^{'}\in A \setminus \{a\} \}=1$.} \vskip 2mm
\noindent \tab L'étude des pgcd des différences d'éléments de
certaines sous-parties d'une base va être un outil important dans
ce qui va suivre. \vskip 5mm \noindent {\bf 2.1.--- Bases
associées à une base additive et dessentialisation élémentaire.}
\vskip 2mm \noindent \tab Supposons donnés une base additive $A$
et $\hbox{\rm Ess}(A) = \{x_1,\cdots,x_s\}$ ($s\in \N$) l'ensemble
de ses éléments essentiels. On appelle {\it diviseur associé} à
l'élément essentiel $x_i$, l'entier positif :
$$d_i=\hbox{\rm pgcd}\left\{x-y ~|~ x,y \in A \setminus \{x_i\}\right\}$$
que l'on sait être $\geq 2$ d'après le lemme 1. Par ailleurs, on sait (voir [DG, lemme 4]) que les entiers $d_1,\cdots,d_s$ sont premiers entre eux deux à deux. Une fois donnés les entiers $d_i$, on pose
$$q(A)=d_1 \cdots d_s$$
(en convenant que $q(A) = 1$ lorsque $A$ est une base sans élément essentiel) et
$$m(A)=\hbox{\rm pgcd}\left\{x-y ~|~ x,y \in A \setminus \hbox{\rm Ess}(A)\right\}$$
\tab Nous appelerons l'entier $m(A)$, le {\it module} de la base $A$ (en particulier $m(A) = 1$ lorsque $A$ est sans élément essentiel). Comme $m(A)$ est clairement multiple de chacun des entiers $d_i$ et que ces derniers sont premiers entre eux deux à deux, l'entier $m(A)$ est donc multiple de $q(A)$. Il est à noter que, de manière générale, $m(A)$ est différent de $q(A)$.
\vskip 2mm
\noindent
\tab Soit $x_0$ le plus petit élément non essentiel de $A$. On appelle {\it ensemble primitif} de $A$, le
sous-ensemble de $\N$ noté $P(A)$ et défini par :
$$P(A)=\left\{\frac{x - x_0}{m(A)} ~|~ x \in A \setminus \hbox{\rm Ess}(A)\right\}$$
et on appelle {\it ensemble élémentaire} associé à $A$, le
sous-ensemble de $\N$, noté $D(A)$ et défini par :
$$D(A)=\left(m(A) \N + x_0\right) \cup \hbox{\rm Ess}(A)$$
\tab Les ensembles $D(A)$ et $P(A)$ sont en fait des bases additives dont on peut estimer l'ordre :
\vskip 2mm
\noindent
{\bf Proposition 2.---} {\it Soient $A$ une base additive, $P(A)$ l'ensemble primitif de $A$ et $D(A)$ l'ensemble élémentaire associé à $A$. Les ensembles $P(A)$ et $D(A)$ sont des bases additives et on a :
$$\hbox{\rm ord}(D(A))\leq \hbox{\rm ord}(A)\leq \hbox{\rm ord}(P(A))+\hbox{\rm ord}(D(A))-1$$}
\vskip 1mm
\noindent
{\bf Preuve :} Comme $D(A)$ contient $A$ il est clair que $D(A)$ est une base. Notons $h$ l'ordre de $A$ et montrons que $P(A)$ est une base additive. Commen\c cons par montrer que $P(A)\cup \{1\}$ est une base additive :
\vskip 2mm
\noindent
\tab Soit $z$ un entier positif assez grand pour que $m(A).z\in h A$. Nous allons établir l'existence d'un entier $h^{'}\in \N^*$, dépendant de $A$ mais pas de $z$, tel que $z$ soit une somme de $h^{'}$ éléments de $P(A) \cup \{1\}$. Par hypothèse, on peut écrire :
$$m(A).z=a_1 x_1 + \cdots + a_s x_s + y_1 + \cdots + y_{\ell}$$
avec $a_1 , \cdots , a_s$ et $\ell$ des entiers positifs vérifiant
$a_1+\cdots+a_s+\ell=h$ et $y_1,\cdots,y_{\ell}$ des éléments de
$A \setminus \{x_1,\cdots,x_s\}$. On a donc :
$$\begin{array}{lll}
z&=&\displaystyle \frac{a_1x_1+\cdots+a_sx_s+y_1+\cdots+y_{\ell}}{m(A)}\\
&=&\displaystyle \frac{y_1-x_0}{m(A)}+\cdots+\frac{y_{\ell}-x_0}{m(A)}+\frac{a_1 x_1+\cdots+a_s x_s+\ell x_0}{m(A)}\\
\end{array}$$
Or, les $\displaystyle
z_1=\frac{y_1-x_0}{m(A)},\cdots,z_{\ell}=\frac{y_{\ell}-x_0}{m(A)}$
sont dans $P(A)$, ce qui implique que le réel
$$u=\frac{a_1x_1+\cdots+a_sx_s+\ell x_0}{m(A)}$$
est un entier positif. Par ailleurs, on a
$$u\leq\frac{h}{m(A)}\max\{x_0,x_1,\cdots ,x_s\}$$
car $a_1+\cdots +a_s+\ell=h$. L'égalité $z=z_1+\cdots +z_{\ell}+u$
est donc une écriture de l'entier $z$ comme somme de $\lfloor
h+\frac{h}{m(A)}\max\{x_0,x_1,\cdots ,x_s\}\rfloor$ éléments de
$P(A)\cup \{1\}$ (car $\ell \leq h$ et $0\in P(A)$). Ainsi,
$P(A)\cup\{1\}$ est une base additive. \vskip 2mm \noindent \tab
Si $1\in P(A)$ alors $P(A)$ est bien une base. Si, maintenant,
$1\notin P(A)$, alors par définition m\^eme du module $m(A)$ on a
$$\hbox{\rm pgcd}\{a-b ~|~ a,b\in P(A)\}=1$$
le lemme 1 assure alors que $P(A)=(P(A)\cup \{1\}) \setminus
\{1\}$ est une base. \vskip 2mm \noindent \tab Montrons maintenant
les encadrements annoncés. On note $H$ l'ordre de $P(A)$ et $\mu$
l'ordre de $D(A)$. \vskip 2mm \noindent $\bullet$ Le fait que
$\hbox{\rm ord}(A) \geq \hbox{\rm ord}(D(A))$ provient de
l'inclusion $A \subset D(A)$. \vskip 2mm \noindent $\bullet$
Montrons maintenant que $\hbox{\rm ord}(A) \leq \hbox{\rm
ord}(P(A))+\hbox{\rm ord}(D(A))-1$. Pour ce faire, nous
choisissons un entier $N_1$ assez grand de fa\c con à ce que :
$$\begin{array}{l}
\bullet N_1 > 2 (H + \mu - 1) \max\{x_0 ,x_1 , \cdots , x_s\}\\
\bullet N_1 - (H - 1) x_0\ \hbox{\rm s'exprime comme une somme de $\mu$ éléments de $D(A)$}\\
\bullet \hbox{\rm tout entier}\ n > \frac{N_1}{2m(A)}\ \hbox{\rm s'exprime comme une somme de $H$ éléments de $P(A)$}
\end{array}$$
\tab Soient $\alpha_1 ,\cdots ,\alpha_s,t,k_1 ,\cdots ,k_t$ des
entiers positifs tels que $\alpha_1+\cdots+\alpha_s+t=\mu$ et tels
que
\begin{equation}\label{e1}
N_1-(H-1)x_0=\alpha_1 x_1+\cdots+\alpha_sx_s+(k_1m(A)+x_0)+\cdots+(k_t m(A)+x_0)
\end{equation}
On a donc :
\begin{eqnarray*}
(k_1 + \cdots k_t)m(A)&=&N_1-\left(\alpha_1 x_1+\cdots +\alpha_s x_s+(t+H-1)x_0\right) \\
&\geq&N_1-(\alpha_1 + \cdots + \alpha_s + t + H - 1) \max\{x_0,\cdots , x_s\}\\
&=&N_1-(H+\mu -1)\max\{x_0,\cdots , x_s\}>\frac{N_1}{2}
\end{eqnarray*}
car $N_1>2(H+\mu -1)\max\{x_0,\cdots ,x_s\}$. Ainsi, on a :
$$k_1 + \cdots + k_t ~>~ \frac{N_1}{2 m(A)}$$
\tab On en déduit que $t\geq 1$ et que le nombre
$(k_1+\cdots+k_t)$ est une somme de $H$ éléments de $P(A)$,
c'est-à-dire que l'on a :
$$k_1+\cdots+k_t=\frac{y_1-x_0}{m(A)}+\cdots+\frac{y_H-x_0}{m(A)}$$
pour certains éléments $y_1 ,\cdots,y_H$ de $A-\{x_1 ,\cdots , x_s\}$. Par suite, la relation (\ref{e1}) donne
\begin{eqnarray*}
N_1 & = & \alpha_1 x_1 + \cdots + \alpha_s x_s + (t + H - 1) x_0 +(k_1 + \cdots + k_t) m(A) \\
& = & \alpha_1 x_1 + \cdots + \alpha_s x_s + (t + H - 1) x_0 + y_1 + \cdots + y_H - H x_0 \\
& = & \alpha_1 x_1 + \cdots + \alpha_s x_s + (t - 1) x_0 + y_1 +\cdots + y_H
\end{eqnarray*}
et montre que $N_1$ est une somme de $\alpha_1 + \cdots + \alpha_s + t - 1 + H=H + \mu - 1$ éléments de $A$ et donc,
$$\hbox{\rm ord}(A) ~\leq~ H + \mu - 1 = \hbox{\rm ord}(P(A)) + \hbox{\rm ord}(D(A)) - 1$$
\fin \newpage \tab On peut aussi encadrer l'ordre de $D(A)$ en
fonction de la base $A$ : \vskip 2mm \noindent {\bf Théorème
3.---} {\it Soient $A$ une base additive, $x_1,\cdots ,x_s$ ses
éléments essentiels, $d_1,\cdots ,d_s$ les diviseurs associés à
$x_1,\cdots ,x_s$, $m(A)$ le module de $A$ et $q(A)=d_1\cdots
d_s$. On a
$$\sum_{i = 1}^{s} d_i - s + 1\leq \hbox{\rm ord}(D(A))\leq \frac{m(A)}{q(A)} \left(\sum_{i = 1}^{s} d_i -
s\sqrt[s]{\frac{q(A)}{m(A)}}\right)+1$$}
\vskip 1mm \noindent {\bf Preuve :} Posons $\mu=\hbox{\rm
ord}(D(A))$ et commen\c cons par établir la minoration annoncée.
Considérons un entier $n$ assez grand tel que l'entier
$$\omega=\sum_{i=1}^{s} (d_i-1)(x_i-x_0)+\mu x_0+n m(A)$$
vérifie $\omega>\mu \max\{x_1,\cdots,x_s\}$ et tel que que
$\omega$ soit somme de $\mu$ éléments de $D(A)$. Il existe donc
des entiers positifs $\alpha_1,\cdots,\alpha_s,k,n_1,\cdots,n_k$
tels que $\alpha_1+\cdots+\alpha_s+k=\mu$ et que
\begin{equation}\label{e2}
\omega=\alpha_1 x_1+\cdots+\alpha_s x_s+(n_1m(A)+x_0)+\cdots+(n_km(A)+x_0)
\end{equation}
On a alors
$$\sum_{i=1}^s (d_i-1)(x_i-x_0)=\sum_{i = 1}^s \alpha_i(x_i-x_0)+(n_1+\cdots+n_k-n)m(A)$$
Pour $i\in \{1,\cdots,s\}$ fixé, en prenant les classes modulo
$d_i$ des deux membres de cette dernière équation, on obtient :
$$(d_i-1) (x_i-x_0)\equiv \alpha_i(x_i-x_0)\ \hbox{\rm mod}(d_i)$$
($x_j-x_0\equiv 0\ \hbox{\rm mod}(d_i)$ pour tout $j \neq i$ et
$d_i|m(A)$). \vskip 2mm \noindent \tab Maintenant $(x_i-x_0)$ est
premier avec $d_i$. En effet, considérons un diviseur $d$ commun à
$x_i-x_0$ et $d_i$. La définition m\^eme de $d_i$ entra\^\i ne que
$d$ divise $x - x_0$ pour tout $x \in A -\{x_i\}$ et comme $d$
divise aussi $x_i - x_0$, alors $d$ divise tous les entiers
$x-x_0$ pour $x$ parcourant $A$. Donc les éléments de $A \setminus
x_0$ sont tous multiple de $d$, mais comme $A \setminus x_0$ est
une base (puisque $A$ en est une et que la basicité est conservée
par translation) cela implique $d=1$. \vskip 2mm \noindent \tab On
en déduit donc que $\alpha_i\equiv (d_i - 1)\ \hbox{\rm mod}(d_i)$
et comme $\alpha_i$ est positif, on a alors
\begin{equation}\label{e3}
\alpha_i\geq d_i-1
\end{equation}
(et ceci étant valable pour tout $1 \leq i \leq s$). \vskip 2mm
\noindent \tab Par hypothèse $\omega - (\alpha_1 x_1 + \cdots +
\alpha_s x_s) > 0$ et donc $k \geq 1$. L'ordre $\mu$ de $D(A)$ est
donc minoré par :
$$\mu=\alpha_1+\cdots +\alpha_s+k\geq \sum_{i=1}^s (d_i-1)+1=\sum_{i = 1}^s d_i-s+1$$
\vskip 2mm
\noindent
\tab Montrons maintenant la majoration. Remarquons préliminairement que
$$\hbox{\rm pgcd}\{m(A),x_1-x_0,\cdots ,x_s-x_0\}=1$$
\tab En effet, si $d \in \N^*$ est un diviseur commun des entiers
$x_i - x_0$ $(1 \leq i \leq s)$ et $m(A)$ alors $d$ divise tout
les entiers $x - x_0$ avec $x \in A$, car si $x \in A \setminus
\hbox{\rm Ess}(A)$, comme $x_0 \in A \setminus \hbox{\rm Ess}(A)$,
le module $m(A)$ de $A$ divise $x - x_0$ et donc a fortiori $d$
divise $x - x_0$. Ainsi, tout élément de l'ensemble $A \setminus
x_0$ est multiple de $d$. Mais comme $A \setminus x_0$ est une
base additive, on a obligatoirement $d = 1$. \vskip 2mm \noindent
\tab Soit $m(A) = q_{1}^{\alpha_1} \cdots q_{k}^{\alpha_k}$ la
décomposition en facteurs premiers de $m(A)$. On note $\Omega$
l'ensemble des diviseurs premiers de $m(A)$, c'est-à-dire $\Omega
= \{q_1,\cdots ,q_k\}$ et $\Omega_i$ $(1\leq i \leq s)$ l'ensemble
des diviseurs premiers de $x_i-x_0$. L'égalité $\hbox{\rm
pgcd}\{m(A),x_1-x_0,\cdots ,x_s-x_0\}=1$ s'écrit en termes
ensemblistes $\Omega \cap \Omega_1 \cap \cdots \cap \Omega_s
=\emptyset$ ou encore
$$\Omega \subset {\complement}_{\mathcal{P}}(\Omega_1 \cap \cdots \cap
\Omega_s)={\complement}_{\mathcal{P}}(\Omega_1) \cup \cdots \cup
{\complement}_{\mathcal{P}}(\Omega_s)$$ o\`u $\cal P$ désigne
l'ensemble des nombres premiers. Ainsi, on a $\Omega=\Gamma_1
\sqcup \cdots \sqcup \Gamma_s$ avec
$$\Gamma_1= \left(\Omega \cap {\complement}_{\mathcal{P}}(\Omega_1)\right),\ \Gamma_i = \left(\Omega \cap
{\complement}_{\mathcal{P}}(\Omega_i)\right)- \bigcup_{j=1}^{i-1} \left(\Omega \cap {\complement}_{\mathcal{P}}(\Omega_j)\right)$$
pour tout $2 \leq i \leq s$. Pour tout $i=1,\cdots ,s$, posons
$$t_i=\prod_{1 \leq r \leq k / q_r \in \Gamma_i} q_{r}^{\alpha_r}$$
\tab Il est clair que les $t_i$ sont des entiers strictement
positifs, deux à deux premiers entre eux et qu'ils vérifient $m(A)
= t_1 \cdots t_s$ et $\hbox{\rm pgcd}\{t_i , x_i-x_0\} = 1$ pour
tout $1 \leq i \leq s$. Montrons que $d_i$ divise $t_i$ pour tout
$i=1,\cdots ,s$. Pour ce faire, nous remarquons que pour tout
couple $(i,j) \in {\{1,\cdots ,s\}}^2$, tel que $i \neq j$, on a
$\hbox{\rm pgcd}(d_i , t_j)=1$ (car $\hbox{\rm pgcd}\{t_j , x_j -
x_0\} = 1$ et $d_i$ divise $(x_j - x_0)$). Ceci entra\^{\i}ne que
pour tout $1\leq i \leq s$, l'entier $d_i$ est premier avec le
produit $t_1 \cdots \widehat{t_i} \cdots t_s$, mais comme $d_i$
divise $(t_1 \cdots \widehat{t_i} \cdots t_s) t_i = m(A)$, on
conclut par le lemme de Gauss que $d_i$ divise $t_i$. \vskip 2mm
\noindent \tab En résumé, on vient de trouver des entiers
$t_1,\cdots ,t_s$ tels que : \vskip 2mm \noindent
$\bullet$ $m(A) = t_1 \dots t_s$\\
$\bullet$ $\hbox{\rm pgcd}\{t_i , t_j\} = 1$ pour tous $1 \leq i < j \leq s$\\
$\bullet$ $\hbox{\rm pgcd}\{t_i , x_i - x_0\} = 1$ pour tout $1 \leq i \leq s$\\
$\bullet$ $d_i$ divise $t_i$ pour tout $1 \leq i \leq s$
\vskip 2mm
\noindent
\tab Considérons alors l'entier
$$\lambda=\sum_{i = 1}^{s} \frac{t_1 \cdots \widehat{t_i} \cdots t_s}{d_1 \cdots \widehat{d_i} \cdots d_s}
(t_i - 1)=\sum_{i = 1}^{s} \frac{m(A)}{q(A)} \left(d_i -
\frac{d_i}{t_i}\right)=\frac{m(A)}{q(A)} \left(\sum_{i = 1}^{s}
d_i - \sum_{i=1}^s \frac{d_i}{t_i}\right)$$ Comme la moyenne
arithmétique de toute famille finie de réels positifs est
supérieure ou égale à sa moyenne géométrique, on a :
$$\frac{1}{s} \sum_{i = 1}^{s} \frac{d_i}{t_i} ~\geq~ \left(\prod_{i = 1}^{s}
\frac{d_i}{t_i}\right)^{\!\!1/s}=\sqrt[s]{\frac{q(A)}{m(A)}}$$
Ceci permet de majorer $\lambda$ de la manière suivante :
$$\lambda \leq \frac{m(A)}{q(A)} \left(\sum_{i = 1}^{s} d_i - s \sqrt[s]{\frac{q(A)}{m(A)}}\right)$$
\tab Nous allons montrer maintenant que tout entier positif assez
grand $N_2$ est une somme d'un nombre $\leq (\lambda + 1)$
éléments de $D(A)$, ce qui prouvera la majoration annoncée dans le
théorème. Soit donc $N_2 \geq (\lambda + 1) \max\{x_0 , x_1 ,
\cdots , x_s\}$ un entier positif arbitraire et $N_3=N_2-(\lambda
+ 1)x_0$. \vskip 2mm \noindent \tab Pour tout $1 \leq i \leq s$,
notons $n_i$ l'unique solution en $n$, dans $\Z \cap [0,t_i-1]$,
de la congruence
$$\left(\frac{x_i - x_0}{d_1 \cdots \widehat{d_i} \cdots d_s} t_1 \cdots \widehat{t_i} \cdots t_s\right)
n \equiv N_3\ \hbox{\rm mod}(t_i)$$ (l'existence et l'unicité de
$n_i$ sont garanties par le fait que $t_i$ est premier avec
$\frac{x_i - x_0}{d_1 \cdots \widehat{d_i} \cdots d_s} t_1 \cdots
\widehat{t_i} \cdots t_s$.) \vskip 2mm \noindent \tab L'entier
$\displaystyle N_3 - \sum_{i = 1}^{s} \left(\frac{x_i - x_0}{d_1
\cdots \widehat{d_i} \cdots d_s} t_1 \cdots \widehat{t_i} \cdots
t_s\right) n_i$ est ainsi multiple de tous les entiers strictement
positifs $t_1 , \cdots , t_s$. Comme ces derniers entiers sont
deux à deux premiers entre eux, ce m\^eme entier est multiple de
leur produit $t_1 \cdots t_s = m(A)$. Ceci entra\^{\i}ne que
l'entier $N_3$ s'écrit sous la forme :
$$N_3=\sum_{i=1}^s \left(\frac{x_i - x_0}{d_1 \cdots \widehat{d_i} \cdots d_s} t_1 \cdots \widehat{t_i}
\cdots t_s\right) n_i + k m(A)$$
pour un certain $k\in \mathbb{Z}$. Ainsi, $N_2$ s'écrit :
\begin{eqnarray*}
N_2&=&\sum_{i = 1}^{s} \left(\frac{x_i-x_0}{d_1 \cdots \widehat{d_i}\cdots d_s} t_1 \cdots
\widehat{t_i}\cdots t_s\right) n_i+k m(A)+(\lambda + 1) x_0 \\
&=&\sum_{i = 1}^{s} \left(\frac{t_1 \cdots \widehat{t_i} \cdots
t_s}{d_1 \cdots \widehat{d_i} \cdots d_s} n_i\right)
x_i+\left(\lambda - \sum_{i = 1}^{s} \frac{t_1 \cdots
\widehat{t_i}\cdots t_s}{d_1 \cdots \widehat{d_i} \cdots d_s}
n_i\right) . x_0 +(k m(A) + x_0)
\end{eqnarray*}
\tab Le fait $N_2 \geq (\lambda + 1) \max\{x_0 , x_1 , \cdots ,
x_s\}$ assure la positivité de l'entier $k$ et montre par suite
que $N_2$ est une somme de $(\lambda + 1)$ éléments de $D(A)$.
\fin \vskip 2mm \noindent
{\bf Remarque :} L'encadrement obtenu
dans ce théorème est en fait optimum. En effet, on voit que
l'encadrement devient une égalité dès que $q(A)=m(A)$ et ceci est
le cas pour le choix de $A=A_n=p_1\cdots p_n\N\cup\{p_1\cdots
\widehat{p_i}\cdots p_n,\ i=1,\cdots ,n\}$. Dans ce cas, on a en
fait $A_n=D(A_n)$ et on en déduit donc que
$$\hbox{\rm ord}A_n=\sum_{i=1}^n p_i-n+1$$
(Ce résultat figure déjà dans [DG] mais y était obtenu de manière
{\it ad hoc}.) \vskip 2mm \noindent \tab Lorsque $A$ est sans
élément essentiel, il est égal à un translaté de son ensemble
primitif $P(A)$. Moralement, l'ensemble primitif associé à une
base $A$ est la base obtenu en "retirant" les éléments essentiels
de $A$. Pourtant, $P(A)$ peut présenter à son tour des éléments
essentiels. On est donc amené à considérer, étant donné une base
$A$, la suite de bases associées $(P^n(A))_n$ définie par
$$P^0(A)=A\ \ \hbox{et}\ \ \forall n\geq 1,\ P^n(A)=P(P^{n - 1}(A))$$
\tab Il est remarquable que, de manière générale, le nombre $s_n$
d'éléments essentiels de $P^n(A)$ ne définit pas forcément une
suite décroissante (voir remarque après la preuve du théorème 4).
Pourtant, la suite $(P^n(A))_n$ est bien stationnaire : \vskip 2mm
\noindent {\bf Théorème 4 .---} {\it Soit $A$ une base additive.
Il existe un entier positif $n$ tel que $P^n(A)$ soit sans élément
essentiel.} \vskip 2mm \noindent {\bf Preuve:} Pour une base
additive donnée $A$, on note $\rho_A$ l'application :
$$\begin{array}{rcl}
\rho_A : P(A) & \longrightarrow & A \\
y & \longmapsto & m(A) . y + x_A\\
\end{array}$$
o\`u $x_A$ désigne le plus petit élément non essentiel de $A$. Il
est clair que l'application $\rho_A$ est injective et que son
image est égale à $A \setminus {\rm Ess}(A)$. \vskip 2mm \noindent
\tab Supposons qu'il existe une base additive $A$ pour laquelle
${\rm Ess}(P^n(A))$ soit non vide pour tout $n \in \mathbb{N}$.
Désignons par $s_i$ $(i \in \mathbb{N})$ le cardinal de l'ensemble
des éléments essentiels de $P^i(A)$. Ces entiers $s_i$ sont donc
tous strictement positifs. Posons, pour tout $n \in \N$,
$\rho_n=\rho_{P^n(A)}$ et considérons la suite infinie
d'applications
$$\cdots \longrightarrow P^{n + 1}(A) \stackrel{\rho_n}{\longrightarrow} P^n(A)
\stackrel{\rho_{n - 1}}{\longrightarrow} P^{n - 1}(A)
\longrightarrow \cdots \longrightarrow P(A)
\stackrel{\rho_0}{\longrightarrow} A$$ ainsi que l'application
composée $f_n$
$$f_n=\rho_0 \circ \rho_1 \circ \dots \circ \rho_n$$
C'est une application affine de $P^{n + 1}(A)$ dans $A$, de pente
$$m(A).m(P(A))\cdots .m(P^n(A))$$
et d'image une partie cofinie de $A$ de cocardinal égal à
$$\#\left({\rm Ess}(A)\right) + \#\left({\rm Ess}(P(A))\right) + \cdots + \#\left({\rm Ess}(P^n(A))\right)=
s_0+s_1+\cdots +s_n$$
\tab Maintenant, comme pour toute base additive $B$, les diviseurs
associés aux éléments essentiels de $B$ sont tous $\geq 2$, on a
$m(B)\geq q(B) \geq 2^{\#({\rm Ess}(B))}$), on en déduit donc que
$$m(A).m(P(A))\cdots m(P^n(A))\geq 2^{s_0}.2^{s_1}\cdots 2^{s_n}= 2^{s_0+s_1+\cdots +s_n}$$
\tab  En posant, pour tout $n \in \N$, $C_n=f_n\left(P^{n +
1}(A)\right) \subset A$ et $u_n = s_0 + s_1 + \dots + s_n$, on a
alors :
$$\left\{
\begin{array}{ll}
\#(A \setminus C_n)=u_n & ~~~~~~~~~~{\rm (I)}\\
{\rm pgcd}\left\{a - b ~|~ a,b \in C_n\right\}\geq 2^{u_n}
&~~~~~~~~~~{\rm (II)}
\end{array}
\right.$$ \tab Pour un entier positif $n$ donné, estimons l'entier
$A(2^{u_n})$ o\`u $A(x)=\# [0,x]\cap A$ pour tout réel $x\geq 0$.
Considérons un entier $a\in A$, tel que $a\leq 2^{u_n}$. Si
$a\notin C_n$ alors, d'après (I), $a$ ne peut prendre que $u_n$
valeurs possibles. Si $a\in C_n$, comme $a\leq 2^{u_n}$ et d'après
la relation (II), $a$ ne peut prendre que deux valeurs possibles.
Ainsi, pour tout $n\in \N$, on a
$$A(2^{u_n})\leq u_n+2$$
\tab Maintenant, un argument combinatoire classique montre que, si
$h$ désigne l'ordre de $A$ et $N$ le plus petit entier positif à
partir duquel tout entier est somme de $h$ éléments de $A$, alors
pour tout $x \geq N$ on a $A(x) \geq (x - N)^{1/h}$. En
particulier, pour tout $\displaystyle n\geq \frac{\log N}{\log
2}$, on a :
$$A(2^{u_n}) \geq (2^{u_n} - N)^{1/h}$$
(car la suite ${(u_n)}_n$ est une suite strictement croissante
d'entiers, elle est donc minorée par $n$ et donc $2^{u_n} \geq 2^n
\geq N$). On en déduit que, pour tout $\displaystyle n\geq
\frac{\log N}{\log 2}$ on a $\displaystyle \frac{\log(u_n +
2)}{\log(2^{u_n} - N)} \geq \frac{1}{h}$. Ce qui est absurde, par
passage à la limite. \fin \vskip 2mm \noindent {\bf Remarque :} a)
Le théorème $4$ explique donc que le calcul de la suite
$(P^n(A))_n$ fournit un algorithme pour ``dessentialiser
élémentairement'' une base, c'est-à-dire un moyen pour associer à
une base une autre base sans élément essentiels. La preuve du
théorème 4, plus exactement l'inégalité  $\displaystyle
\frac{\log(u_n + 2)}{\log(2^{u_n} - N)} \geq \frac{1}{h}$, permet
de majorer le nombre $\delta(A)$ d'étapes nécessaire à la
dessentialisation élémentaire de $A$ en fonction des invariants
$h$ et $N$ liés à $A$. Par exemple, une simple étude de fonction
montre que
$$\delta(A) \leq \max\left(\frac{\log N}{\log 2}+1,N^{1/h}-2\right)$$
b) Une conséquence immédiate du théorème 4 est que, si $A$ est une base additive possédant des éléments essentiels, il existe une base $B$ sans élément essentiel et deux entiers $a,b$ avec $a\geq 2$ tels que
$$A\sim aB+b$$
Ce résultat sera généralisé au paragraphe 3 dans le cas des parties essentielles.
\vskip 2mm
\noindent
c) Il est intéressant de remarquer qu'étant donné un entier $n\geq 0$ et une suite d'entiers strictement positifs $s_0,\cdots ,s_n$ il existe toujours une base additive $A$ telle que $\delta(A)=n$ et tel que pour tout $i=0,\cdots ,n$, le nombre d'éléments essentiels de $P^i(A)$ soit égal à $s_i$. En effet, considérons la suite finie de bases additives, définie par récurrence de la manière suivante :
$$A_0=\N,\ \forall i=0,\cdots n,\ A_{i+1}=p_1\cdots p_{s_{n-i}}A_i\cup \left\{p_1\cdots \widehat{p_j}\cdots p_{s_{n-i}},\ j=1,\cdots,s_{n-i} \right\}$$
\tab Les ensembles $A_i$ sont bien des bases et elles comptent respectivement $s_{n-i+1}$ éléments essentiels. Il est clair que l'on a, pour tout $i=0,\cdots ,n$, $P(A_{i+1})=A_i$
\vskip 5mm
\noindent
{\bf 2.2.--- Etude du nombre d'éléments essentiels.}
\vskip 2mm
\noindent
\tab Une conséquence intéressante de la partie précédente est le lemme suivant, qui permet de relier l'ordre d'une base à son nombre d'éléments essentiels :
\vskip 2mm
\noindent
{\bf Lemme 6.---} {\it Soit $A$ une base additive d'ordre $h$, possédant $s$ éléments essentiels. On a
$$h\geq \sum_{i=1}^s p_i-s+1$$}
\vskip 1mm \noindent {\bf Preuve:} D'après la proposition 2 et le
théorème 3, on a $\hbox{\rm ord}(A)\geq \hbox{\rm ord}(D(A))\geq
\sum_{i=1}^s d_i-s+1$ et comme les $d_i$ sont premiers entre eux
deux à deux et supérieurs à $2$, on a $\sum_{i=1}^s d_i\geq
\sum_{i=1}^s p_i$. \fin \vskip 2mm \noindent \tab Si l'on
considère la base $A_n=p_1\cdots p_n\N\cup\{p_1\cdots
\widehat{p_i}\cdots p_n,\ i=1,\cdots ,n\}$ on sait que son ordre
est $h_n=\sum_{i=1}^{n}p_i-n+1$ et qu'elle compte exactement $n$
éléments essentiels. On voit donc que, pour le choix $A=A_n$,
l'inégalité du Lemme 6 est optimale. \vskip 2mm \noindent \tab Le
résultat principal de [DG] assure que si $A$ est une base additive
d'ordre $h$ comptant $s$ éléments essentiels, alors on a
l'inégalité
$$s\leq 5,7\sqrt{\frac{h}{\log h}}$$
De plus, il est prouvé dans [DG] qu'il existe une constante $c>0$, indépendante de $n$, telle que
$$n\geq c\sqrt{\frac{h_n}{\log h_n}}$$
Ainsi, on savait depuis [DG] que le comportement de $s$ en fonction de $h$ était en $O(\sqrt{\frac{h}{\log h}})$ et qu'il s'agissait l\`a de la meilleure estimation possible en $O$. Il restait deux questions en suspend : la détermination de la constante absolu, c'est-à-dire le plus petit réel $C>0$ tel que pour toute base additive d'ordre $h$ comptant $s$ éléments essentiels on ait $s\leq C\sqrt{\frac{h}{\log h}}$ et la détermination de la constante asymptotique, c'est-à-dire le réel (s'il existe) $c=\lim_h s(h)\sqrt{\frac{\log h}{h}}$ o\`u $s(h)$ désigne le nombre maximal d'éléments essentiels d'une base additive d'ordre $h$. Si l'on regarde la suite $(A_n)_n$, on voit qu'asymptotiquement on a $n\sqrt{\frac{\log h_n}{h_n}}\simeq 2$, ce qui suggère que la constante asymptotique au problème est $2$. A.Plagne a montré récemment que c'était bien le cas (cf [Pl]). Le théorème suivant fournit la constante absolu et montre comment obtenir la constante asymptotique supérieure :
\vskip 2mm
\noindent
{\bf Théorème 7.---} {\it a) On a $\displaystyle \limsup_h s(h)\sqrt{\frac{\log h}{h}}=2$.
\vskip 2mm
\noindent
b) Soit $A$ une base additive d'ordre $h\geq 2$ comptant $s$ éléments essentiels.
\vskip 2mm
\noindent
\tab $\bullet$ On a
$$s\leq C\sqrt{\frac{h}{\log h}}$$
avec $C=30\sqrt{\frac{\log 1564}{1564}}\sim 2,0572841285$. De plus, cette inégalité devient une égalité pour $A=A_{30}=p_1\cdots p_{30}\N\cup\{p_1\cdots \widehat{p_i}\cdots p_{30},\ i=1,\cdots ,30\}$.
\vskip 2mm
\noindent
\tab $\bullet$ Pour tout $\alpha >2$ on a
$$s\leq \alpha \sqrt{\frac{h}{\log h}}$$
dès que $h\geq \exp\left(-\frac{7}{2}\alpha^2\frac{\log(\alpha^2-4)}{(\alpha^2-4)}\right)$.}
\vskip 2mm
\noindent
{\bf Preuve:} a) Avec les notations précédentes, notons $C_h=s(h)\sqrt{\frac{\log h}{h}}$. Quand $h=h_n$, on sait que la base additive $A_n$, d'ordre $h_n$, compte exactement $n$ éléments essentiels, on a donc $s(h_n)\geq n$. Maintenant, par le lemme 6, on a
$$\sum_{i=1}^{n}p_i-n+1=h_n\geq \sum_{i=1}^{s(h_n)}p_i-s(h_n)+1$$
ce qui montre que $s(h_n)\leq n$ et, par suite, que $s(h_n)=n$. On a donc $n\sqrt{\frac{\log h_n}{h_n}}=C_{h_n}$.
\vskip 2mm
\noindent
\tab Soit $A$ une base additive d'ordre $h\geq 4$, $n\geq 2$ l'unique entier tel que $h_n\leq h<h_{n+1}$ et $s$ le nombre d'éléments essentiels de $A$. On a
$$\sum_{i=1}^{n+1}p_i-(n+1)+1=h_{n+1}>h\geq \sum_{i=1}^{s}p_i-s+1$$
on en déduit donc que $s<n+1$. Ainsi on a
$$s\sqrt{\frac{\log h}{h}}\leq n\sqrt{\frac{\log h}{h}}\leq n\sqrt{\frac{\log h_n}{h_n}}=C_{h_n}$$
et, par suite, $C_h\leq C_{h_n}$.\\
\tab Il s'ensuit que $\displaystyle \limsup_h s(h)\sqrt{\frac{\log h}{h}}=\limsup_n C_{h_n}= \limsup_n n\sqrt{\frac{\log h_n}{h_n}}$. Le théorème des nombres premiers a comme conséquence l'équivalence suivante : $p_n\simeq_n n\log n$. On en déduit donc que
$$h_n=\sum_{i=1}^np_i-n+1\simeq_n \sum_{i=1}^np_i\simeq_n \sum_{i=1}^ni\log i\simeq_n \frac{1}{2}n^2\log n$$
et, par suite,
$$n\sqrt{\frac{\log h_n}{h_n}}\simeq_n n\sqrt{\frac{2\log n}{\frac{1}{2}n^2\log n}}=2$$
ce qui prouve le a).
\vskip 2mm
\noindent
b) $\bullet$ On cherche donc à déterminer $\max_h C_h$. Dans la preuve du a) on a vu que pour tout entier $h\geq 4$ il existait $n\geq 2$ tel que $C_h\leq C_{h_n}$. Comme par ailleurs, $C_2,C_3\leq C_4$, il s'ensuit que $\max_h C_h=\max_n C_{h_n}$. On va montrer que
$$\max_n C_{h_n}=C_{h_{30}}=C$$
ce qui prouvera l'assertion. \vskip 2mm \noindent \tab Si $n\leq
127042$ on vérifie que $C_{h_n}\leq C$. (Nous ne détaillons pas
ici le calcul qui s'obtient tr\`es facilement avec une table de
nombres premiers et un tableur.) \vskip 2mm \noindent \tab
Supposons $n\geq 127042$. D'après [MR], on sait que l'on a
$$\sum_{i=1}^np_i\geq \frac{n^2}{2}(\log n+\log\log n-1,5034)$$
Ainsi, pour montrer que $n\sqrt{\frac{\log h_n}{h_n}}\leq C$, il suffit de prouver que
$$\frac{n^2}{C^2}\leq \frac{\frac{n^2}{2}(\log n+\log\log n-1,5034)-n+1}{\log
(\frac{n^2}{2}(\log n+\log\log n-1,5034)-n+1)}$$
ce qui équivaut, une fois les simplifications effectuées, à l'inéquation
\begin{center}
\flushleft{$\displaystyle (C^2-4)\log n+(C^2-2)\log\log n-C^2.1,5034-\frac{2C^2}{n}+\frac{2C^2}{n^2}$}\\
\flushright{$\displaystyle -2\log\left(\frac{1}{2}+\frac{\log\log n}{2\log n}-\frac{1,5034}{2\log n}
-\frac{1}{n\log n}+\frac{1}{n^2\log n}\right)\geq 0$}
\end{center}
Comme $n\geq 127042$, on a $\displaystyle \frac{\log\log n}{2\log n}-\frac{1,5034}{2\log n}
-\frac{1}{n\log n}+\frac{1}{n^2\log n}\leq \frac{1}{2}$. Il suffit donc de vérifier que
$$(C^2-4)\log n+(C^2-2)\log\log n-C^2.1,5034-\frac{2C^2}{n}\geq 0$$
pour obtenir le résultat. Compte-tenu du fait que la fonction
incriminée dans le membre de gauche est croissante et que
l'inégalité est vraie pour $n=127042$, on en déduit qu'elle est
vraie pour tout $n\geq 127042$. \vskip 2mm \noindent $\bullet$
Soit $\alpha>2$ (on suppose $\alpha\leq C$ sinon la proposition
est évidente) et $A$ une base additive d'ordre $h$ comptant $s$
éléments essentiels tels que
$$s>\alpha\sqrt{\frac{h}{\log h}}$$
On a $h\geq \sum_{i=1}^sp_i-s+1$. D'après [MR], on sait que
$\sum_{i=1}^sp_i\geq \frac{1}{2}s^2\log s$ pour tout $s\geq 2$, et
donc
$$\begin{array}{ll}
&\displaystyle h> \frac{1}{2}\alpha^2\frac{h}{\log h}\left( \frac{1}{2}\log h-\frac{1}{2}\log\log h+\log
\alpha\right)-\alpha\sqrt{\frac{h}{\log h}}+1\\
\Longleftrightarrow &\displaystyle (4-\alpha^2)\log h+\alpha^2\log\log h-2\alpha^2\log \alpha+4\alpha
\sqrt{\frac{\log h}{h}}-4\frac{\log h}{h}>0\\
\Longrightarrow&\displaystyle (4-\alpha^2)\log h+\alpha^2\log\log h-2\alpha^2\log \alpha+4\alpha>0\\
\end{array}$$
Considérons la fonction $f(x)=(4-\alpha^2)x+\alpha^2\log
x-2\alpha^2\log \alpha+4\alpha$. Elle présente un maximum absolu
en $x_\alpha=\frac{\alpha^2}{\alpha^2-4}$ et est décroissante sur
$[x_\alpha,+\infty[$. On a
$$\begin{array}{lll}
\displaystyle f\left(-\frac{7}{2}\alpha^2\frac{\log(\alpha^2-4)}{(\alpha^2-4)}\right)&=&\displaystyle
\frac{5}{2}\alpha^2\log(\alpha^2-4)+\alpha^2\log(-\log(\alpha^2-4))+\alpha^2\log \frac{7}{2}+4\alpha\\
&<&\displaystyle \alpha^{2}\left(\frac{5}{2}\log(\alpha^2-4)+\log(-\log(\alpha^2-4))+\log \frac{7}{2}+2\right)\\
&=&\displaystyle \alpha^{2}\log\left(-\frac{7e^2}{2}t^{5/2}\log t\right)\ \hbox{\rm avec}\ t=\alpha^2-4\\
\end{array}$$
Une étude rapide de la fonction $\displaystyle t\mapsto
t^{5/2}\log t$ montre, qu'avec l'hypothèse $\alpha \leq C$, on a
alors
$$f\left(-\frac{7}{2}\alpha^2\frac{\log(\alpha^2-4)}{(\alpha^2-4)}\right)<0$$
et donc que, puisque $f(\log h)>0$, on a $\displaystyle \log
h<-\frac{7}{2}\alpha^2\frac{\log(\alpha^2-4)}{(\alpha^2-4)}$
c'est-à-dire $\displaystyle
h<\exp\left(-\frac{7}{2}\alpha^2\frac{\log(\alpha^2-4)}{(\alpha^2-4)}\right)$.
Ceci achève la preuve du théorème. \fin \vskip 10mm \noindent
\centerline{\bf 3.--- Sur les parties essentielles d'une base
additive.} \vskip 2mm \noindent \tab On s'intéresse dans ce
paragraphe aux parties essentielles d'une base additive qui, nous
le rappelons, sont les parties \underline{finies} et
\underline{minimales} $P\subset A$ telles que $A \setminus P$ ne
soit plus une base. Le point de départ pour leur étude est le lien
qu'il y a entre l'existence d'une partie essentielle pour une base
et le fait d'être en progression arithmétique : \vskip 2mm
\noindent {\bf Proposition 8.---} {\it Soit $A$ une base. Les
propriétés suivantes sont équivalentes : \vskip 2mm \noindent i)
$A$ possède une partie essentielle, \vskip 2mm \noindent ii) Il
existe deux entiers $a\geq 2$, $b\geq 0$ et une partie $X\subset
\N$ telle que $A\sim aX+b$.} \vskip 2mm \noindent {\bf Preuve :}
$i)\Rightarrow ii)$ Soit $P$ une partie essentielle de $A$ et
$x\in P$. L'ensemble $(A \setminus P)\cup \{x\}$ est une base
additive possédant $x$ pour élément essentiel, donc d'après le
lemme 1 on a ii). \vskip 2mm \noindent $ii)\Rightarrow i)$ Si
$A\sim aX+b$, on peut écrire $A=(aB+b)\cup E$ avec $E$ finie non
vide (sinon $A$ n'est pas une base). L'ensemble $A \setminus E$
n'est pas une base, on en déduit donc qu'il existe une partie
essentielle $P$ de $A$ telle que $P\subset E$. \fin \vskip 2mm
\noindent \tab En d'autres termes, l'existence de parties
essentielles pour une base équivaut pour cette base à vivre dans
une progression arithmétique non triviale à partir d'un certain
rang. On va voir ici que la raison d'une telle progression donne
des informations capitales pour l'étude des parties essentielles.
\vskip 2mm \noindent \tab Si $A\sim aX+b$ et que l'on impose la
condition $0\leq b<a$ (ce qui est toujours possible quitte à
translater $X$) on voit qu'en posant
$$\begin{array}{l}
\displaystyle E=\left\{x\in A ~|~ x\not\equiv b\ \hbox{\rm mod}(a)\right\}\\
\displaystyle B=\left\{\frac{x-b}{a} ~|~ x\in A,\ x\equiv b\ \hbox{\rm mod}(a)\right\}\\
\end{array}$$
on a $A=(aB+b)\cup E$. On remarque que, toujours sous la condition
$0\leq b<a$, les ensembles $B$ et $E$ sont uniques dés que l'on
impose la condition $x\not \equiv b\ \hbox{\rm mod}(a)$ pour tout
$x\in E$, et que dans cette situation on a $A=(aB+b)\sqcup E$ avec
$E$ fini. \vskip 2mm \noindent \tab \underline{Dans la suite de ce
texte}, quand on supposera que $A\sim aX+b$ et que l'on écrira
alors $A=(aB+b)\cup E$ on supposera toujours que $0\leq b<a$ et
que tout élément $x\in E$ vérifie $x\not \equiv b\ \hbox{\rm
mod}(a)$. \vskip 2mm \noindent {\bf Théorème 9.---} {\it Soit $A$
une partie de $\N$. S'il existe une infinité d'entiers $a\in \N$
tels que $A\sim aX+b$ (pour un certain $b\in \N$ et un certain
$X\subset \N$ dépendants de $a$), alors $A$ n'est pas une base
additive.} \vskip 2mm \noindent {\bf Preuve :} Commen\c cons par
remarquer que s'il existe $a\ne a^{'}$ tels que $A\sim aX+b\sim
a^{'}X^{'}+b^{'}$ alors $A\sim mX^{''}+b^{''}$ o\`u $m=\hbox{\rm
ppcm} (a,a^{'})$ pour une certaine partie $X^{''}$ et un certain
entier $b^{''}$. En particulier, s'il existe une infinité
d'entiers $a$ tels que $A\sim aX+b$ alors il en existe une
infinité ordonnée pour la relation de divisibilité. \vskip 2mm
\noindent \tab Si $a_1,a_2,\cdots $ désigne une suite d'entiers
strictement croissante (pour la relation de divisibilité) tels que
$A\sim a_iX_i+b_i$ alors, avec notre convention, si l'on écrit
$A=(a_iB_i+b_i)\cup E_i$, on a pour tout $i$
$$\left\{\begin{array}{l}
\displaystyle (a_{i+1}B_{i+1}+b_{i+1})\subset (a_iB_i+b_i)\\
\displaystyle E_i\subset E_{i+1}\\
\end{array}
\right.$$ Dans cette situation, l'intersection $\displaystyle
\bigcap_{i}(a_iB_i+b_i)$ est réduite au plus à un seul élément, ce
qui entraine que $A=\bigcup_i E_i\cup \{\alpha\}$ pour un certain
élément $\alpha\in A$. Ainsi, puisque $A$ est infini, la suite
croissante d'ensembles finis $(E_i)_i$ ne peut pas être
stationnaire. On déduit donc que si $A\sim aX+b$ pour une infinité
de $a$, alors pour tout entier $a$ donné tel que $A=(aB+b)\cup E$
et tout entier $M$, il existe un entier $a^{'}>M$ tel que
$a|a^{'}$, $A=(a^{'}B^{'}+b^{'})\cup E^{'}$ et $E\varsubsetneq
E^{'}$ (et donc $(a^{'}B^{'}+b^{'})\varsubsetneq (aB+b)$). \vskip
2mm \noindent \tab Raisonnons par l'absurde en supposant que $A$
soit une base additive. \vskip 2mm \noindent $\bullet$ Soient
$a_1,b_1,B_1,E_1$ tels que $A=(a_1B_1+b_1)\cup E_1$ et $a_1\geq
2$. Comme $A$ est une base, on a $E_1\ne \emptyset$. On considère
un élément $x_1\in E_1$. \vskip 2mm \noindent $\bullet$ Il existe
$a_2,b_2,B_2,E_2$ tels que
$$\left\{
\begin{array}{l}
a_1|a_2\\
a_2>\max E_1\\
A=(a_2B_2+b_2)\cup E_2\\
\end{array}
\right.$$
Il existe $a_3,b_3,B_3,E_3$ tels que
$$\left\{
\begin{array}{l}
a_2|a_3\\
A=(a_3B_3+b_3)\cup E_3\\
E_2\varsubsetneq E_3\\
\end{array}
\right.$$
On peut alors choisir un élément $x_2\in (a_2B_2+b_2)$ tel que $x_2\in E_3$.
\vskip 2mm
\noindent
$\bullet$ De même, il existe $a_4,b_4,B_4,E_4$ tels que
$$\left\{
\begin{array}{l}
a_3|a_4\\
a_4>2\max E_3\\
A=(a_4B_4+b_4)\cup E_4\\
\end{array}
\right.$$
Il existe alors $a_5,b_5,B_5,E_5$ tels que
$$\left\{
\begin{array}{l}
a_4|a_5\\
A=(a_5B_5+b_5)\cup E_5\\
E_4\varsubsetneq E_5\\
\end{array}
\right.$$ On peut choisir un élément $x_3\in (a_4B_4+b_4)$ tel que
$x_3\in E_5$. \vskip 2mm \noindent $\bullet$ Par récurrence, si
l'on suppose au rang $k-1\geq 1$ avoir construit
$a_{2k-2},a_{2k-1}$, $b_{2k-2},b_{2k-1}$, $B_{2k-2},B_{2k-1}$,
$E_{2k-2},E_{2k-1}$ et $x_k$ tels que
$$\left\{\begin{array}{l}
a_{2k-2}>(k-1)\max E_{2k-3}\\
a_{2k-2}|a_{2k-1}\\
A=(a_{2k-2}B_{2k-2}+b_{2k-2})\cup E_{2k-2}=(a_{2k-1}B_{2k-1}+b_{2k-1})\cup E_{2k-1}\\
x_k\in (a_{2k-2}B_{2k-2}+b_{2k-2})\cap E_{2k-1}\\
\end{array}\right.$$
alors au rang $k$ il existe $a_{2k},b_{2k},B_{2k},E_{2k}$ tels que
$$\left\{
\begin{array}{l}
a_{2k-1}|a_{2k}\\
a_{2k}>k\max E_{2k-1}\\
A=(a_{2k}B_{2k}+b_{2k})\cup E_{2k}\\
\end{array}
\right.$$
Il existe alors  $a_{2k+1},b_{2k+1},B_{2k+1},E_{2k+1}$ tels que
$$\left\{
\begin{array}{l}
a_{2k}|a_{2k+1}\\
A=(a_{2k+1}B_{2k+1}+b_{2k+1})\cup E_{2k+1}\\
E_{2k}\varsubsetneq E_{2k+1}\\
\end{array}
\right.$$
On peut alors choisir un élement $x_{k+1}\in (a_{2k}B_{2k}+b_{2k})$ tel que $x_{k+1}\in E_{2k+1}$.
\vskip 2mm
\noindent
\tab Supposons que $h$ soit l'ordre de $A$ et considérons l'entier
$$x=x_1+\cdots +x_h+\lambda a_{2h}$$
o\`u $\lambda$ est un paramètre entier. Puisque $A$ est une base
d'ordre $h$, pour tout $\lambda$ assez grand il existe $y_1,\cdots
,y_h\in A$ tel que
$$x=y_1+\cdots +y_h$$
\tab Nous allons montrer, par récurrence, que l'on peut réordonner
la famille $(y_i)_i$ de sorte que $y_i\in E_{2i-1}$ pour tout
$i=1,\cdots ,h$. \vskip 2mm \noindent $\bullet$ Par construction,
on a $x_p\equiv b_1\ \hbox{\rm mod}(a_1)$ pour tout $p=2,\cdots
,h$ et $x_1\not \equiv b_1\ \hbox{\rm mod} (a_1)$. Ainsi, puisque
$\lambda a_{2h}\equiv 0\ \hbox{\rm mod} (a_1)$ on a
$$y_1+\cdots +y_h=x=x_1+\cdots +x_h+\lambda a_{2h}\not \equiv hb_1\ \hbox{\rm mod} (a_1)$$
et donc il existe un $y_k$, disons $y_1$, tel que $y_1\not \equiv
b_1\ \hbox{\rm mod} (a_1)$, c'est-à-dire $y_1\in E_1$. La
propriété est donc vérifiée au rang $l=1$. \vskip 2mm \noindent
$\bullet$ Supposons que pour $l<h$ on ait $y_i\in E_{2i-1}$ pour
tout $i=1,\cdots ,l$. En particulier, on a $y_1,\cdots ,y_l\in
E_{2l-1}$. Alors au rang $l+1$, on a :
$$y_{l+1}+\cdots +y_h=(x_1-y_1)+\cdots +(x_l-y_l)+x_{l+1}+\cdots +x_h+\lambda a_{2h}$$
Posons
$$\omega=(x_1-y_1)+\cdots +(x_l-y_l)$$
\tab Si $\omega=0$ alors
$$y_{l+1}+\cdots +y_h=x_{l+1}+\cdots +x_h+\lambda a_{2h}$$
Par hypothèse, $x_{l+1}\not\equiv b_{2l+1}\ \hbox{\rm
mod}(a_{2l+1})$ et comme on a $x_p\equiv b_{2l+1}\ \hbox{\rm
mod}(a_{2l+1})$ pour tout $p=l+2,\cdots ,h$ et $\lambda
a_{2h}\equiv 0\ \hbox{\rm mod}(a_{2l+1})$, on en déduit que
$$y_{l+1}+\cdots +y_h\not\equiv (h-l)b_{2l+1}\ \hbox{\rm mod}(a_{2l+1})$$
donc il existe un $y_k$, disons $y_{l+1}$ tel que
$y_{l+1}\not\equiv b_{2l+1}\ \hbox{\rm mod}(a_{2l+1})$,
c'est-à-dire $y_{l+1}\in E_{2l+1}$. \vskip 2mm \noindent \tab Si
$\omega\ne 0$ alors comme $x_1,\cdots ,x_l,y_1,\cdots ,y_l\in
E_{2l-1}$ et que $a_{2l}>l\max E_{2l-1}$, on a $a_{2l}>|
(x_1-y_1)+\cdots +(x_l-y_l)|$. Il s'ensuit que $(x_1-y_1)+\cdots
+(x_l-y_l)\not\equiv 0\ \hbox{\rm mod}(a_{2l})$. Par ailleurs,
pour tout $p=l+1,\cdots ,h$ on a $x_p\equiv b_{2l}\ \hbox{\rm
mod}(a_{2l})$ et $\lambda a_{2h}\equiv 0\ \hbox{\rm mod}(a_{2l})$,
ainsi
$$y_{l+1}+\cdots +y_h\not\equiv (h-l) b_{2l}\ \hbox{\rm mod}(a_{2l})$$
et donc il existe un $y_k$, disons $y_{l+1}$ tel que
$y_{l+1}\not\equiv b_{2l}\ \hbox{\rm mod}(a_{2l})$, c'est-à-dire
$y_{l+1}\in E_{2l}\subset E_{2l+1}$. \vskip 2mm \noindent \tab La
récurrence est ainsi achevée. Nous venons donc de montrer que
l'entier $x$ s'écrit comme somme de $h$ éléments de $E_{2h+1}$ et
ceci indépendamment du paramètre $\lambda$. Comme $E_{2h+1}$ est
un ensemble fini, il n'y a qu'un nombre fini d'entiers $x$ pouvant
s'écrire comme somme de $h$ éléments de $E_{2h+1}$. On en déduit
une absurdité puisque $\lambda$ peut prendre une infinité de
valeurs. Ainsi $A$ n'est pas une base additive. \fin \vskip 2mm
\noindent \tab On déduit du théorème 9 que pour toute base $A$ il
existe un plus grand entier $a$ tel que $A\sim aX+b$. Une fois
donné ce plus grand entier $a$, si $a\geq 2$ il existe alors un
unique entier $0\leq b<a$ et deux parties uniques $B,E$ tels que
$$A=(aB+b)\cup E$$
avec pour tout $x\in E$, $x\not \equiv b\ \hbox{\rm mod}(a)$.
\vskip 2mm \noindent {\bf Définition.---} {\it Avec les notations
précédentes, on appelle l'entier $a$ {\rm la raison} de $A$, la
partie $B$ {\rm la dessentialisée} de $A$ et la partie $E$ {\rm le
réservoir} de $A$.} \vskip 2mm \noindent \tab Par exemple
l'ensemble, $P(k)=\{n^k ~|~ n\in \N\}$ est une base de raison $1$
et de réservoir vide. L'ensemble $\cal P$ des nombres premiers est
une base de raison $2$, son réservoir est $\{2\}$. \vskip 2mm
\noindent \tab Il est à noter qu'il n'existe pas {\it a priori} de
lien entre l'ordre et la raison d'une base. Par exemple, si l'on
considère la suite de bases additives
$\left(aP(k)\cup\{1\}\right)_k$, on voit que la raison de ces
bases est toujours $a$ et que l'ordre tend vers $+\infty$.
Réciproquement, si l'on considère la suite de bases additives
$\left(k\N\cup\{1,\cdots ,k-1\}\right)_k$, on voit que l'ordre de
ces bases est toujours $2$ alors que la raison tend vers
$+\infty$. \vskip 2mm \noindent \tab La proposition suivante
généralise la remarque du théorème 4 : \vskip 2mm \noindent {\bf
Proposition 10.---} {\it a) La dessentialisée d'une base additive
est une base additive sans partie essentielle. \vskip 2mm
\noindent b) Soit $A$ une base additive de raison $a$ et $B$ la
dessentialisée de $A$. Soit $a^{'},b^{'}\in \N$ et $B^{'}\subset
\N$ tels que $A\sim a^{'}B^{'}+b^{'}$. Les propriétés suivantes
sont équivalentes : \vskip 2mm \noindent i) $a=a^{'}$, \vskip 2mm
\noindent ii) $B\sim (B^{'}+k)$, pour un certain entier $k$,
\vskip 2mm \noindent iii) $B^{'}$ est une base sans partie
essentielle.} \vskip 2mm \noindent {\bf Preuve :} a)  Soit $h$
l'ordre de $A$. Pour tout entier $n$ assez grand il existe
$x_1,\cdots ,x_k\in E$, $m_1,\cdots ,m_l\in B$ avec $k+l=h$ tel
que
$$an+b=x_1+\cdots +x_k+(am_1+b)+\cdots +(am_l+b)$$
on en déduit que
$$n=m_1+\cdots+m_l+\frac{x_1+\cdots +x_k+(l-1)b}{a}$$
Maintenant
$$\frac{x_1+\cdots +x_k+(l-1)b}{a}\leq h\frac{\max E+b}{a}$$
on en déduit donc que $B\cup \{0,1\}$ est une base additive
d'ordre $\leq h+z$ o\`u $z$ désigne la partie entière de
$h\frac{\max E+b}{a}$. \vskip 2mm \noindent \tab Supposons que $B$
ne soit pas une base, alors $B\cup\{0,1\}$ possède une partie
essentielle et donc il existe $a^{'}\geq 2$, $b^{'}\in \N$ et
$B^{'}\subset \N$ tels que $B\cup\{0,1\}\sim a^{'}B^{'}+b^{'}$.
Par suite, on a $B\sim a^{'}B^{'}+b^{'}$ et donc $A\sim
a^{'}aB^{'}+c$ pour un certain $c$, ce qui nie la maximalité de
$a$. Donc $B$ est bien une base. \vskip 2mm \noindent \tab Si $B$
possède une partie essentielle alors  $B\sim a^{'}B^{'}+b^{'}$
avec $a^{'}\geq 2$ et donc $A\sim a^{'}aB^{'}+c$, ce qui nie à
nouveau la maximalité de $a$. \vskip 2mm \noindent b)
$i)\Rightarrow ii)$ Il existe un entier $n_0$ tel que
$$(aB+b)\cap [n_0,+\infty[=A\cap [n_0,+\infty[=(aB^{'}+b^{'})\cap [n_0,+\infty[$$
on en déduit que $b\equiv b^{'}\ \hbox{\rm mod}(a)$. Posons $\displaystyle k=\frac{b^{'}-b}{a}$, on a
$$B\cap \left[\frac{n_0-b}{a},+\infty \right[=(B^{'}+k)\cap \left[ \frac{n_0-b}{a},+\infty \right[$$
c'est-à-dire $B\sim B^{'}+k$. \vskip 2mm \noindent $ii)\Rightarrow
iii)$ Soit $n_0$ est tel que $B\cap [n_0,+\infty[=(B^{'}+k)\cap
[n_0,+\infty[$. Comme $B$ est sans partie essentielle, $B\cap
[n_0,+\infty[$ est une base, donc $(B^{'}+k)\cap [n_0,+\infty[$ et
$(B^{'}+k)$ le sont aussi. Si $P$ était une partie essentielle de
$B^{'}+k$ alors $(B^{'}+k)-P$ ne serait pas une base, idem pour
$(B^{'}+k)\cap [n_0,+\infty[-P=B\cap [n_0,+\infty[-P$, et par
suite $[0,n_0[\cup P$ contiendrait une partie essentielle de $B$
ce qui est absurde. Ainsi, $(B^{'}+k)$ est une base sans partie
essentielle, il en est donc de même de sa translatée
$B^{'}=(B^{'}+k)-k$. \vskip 2mm \noindent $iii)\Rightarrow i)$
Supposons $a\ne a^{'}$. Par maximalité de $a$ on a $a^{'}|a$
strictement. Posons $a=ka^{'}$ avec $k\geq 2$, on a donc
$B^{'}\sim kB+t$ pour un certain entier $t$ et par suite
(proposition 8), $B^{'}$ possède une partie essentielle, ce qui
est absurde. \fin \vskip 2mm \noindent \tab Le théorème suivant
montre que toute partie essentielle d'une base $A$ vit dans son
réservoir, ce qui explique la terminologie que nous employons pour
désigner cet invariant. \vskip 2mm \noindent {\bf Théorème 11.---}
{\it Soit $A$ une base additive de raison $a$ et de réservoir $E$.
Si $P$ est une partie essentielle de $A$ alors $P\subset E$. En
particulier, l'ensemble des parties essentielles de $A$ est fini,
son cardinal est majoré par la longueur du radical de $a$ (i.e. le
nombre de nombres premiers divisant $a$).} \vskip 2mm \noindent
{\bf Preuve :} On a $A=(aB+b)\cup E$. Soit $P$ une partie
essentielle de $A$ telle que $P\not \subset E$ et $x_0\in P\cap
(aB+b)$. L'ensemble $(A \subset P)\cup \{x_0\}$ est donc une base
additive possédant $x_0$ pour élément essentiel. Le lemme 1,
permet alors d'affirmer que
$$\hbox{\rm pgcd}\{x-y ~|~ x,y\in A \setminus P\}=d\geq 2$$
\tab Montrons par l'absurde que $d$ ne divise pas $a$. Supposons
donc le contraire et prenons deux éléments $x,y\in (A \setminus
P)\cup \{x_0\}$. \vskip 2mm \noindent
$\bullet$ Si $x,y\in aB+b$ alors $d|x-y$ par hypothèse.\\
$\bullet$ Si $x,y\in (A \setminus P)$ on a $d|(x-y)$ par définition de $d$.\\
$\bullet$ Si $x=x_0$ et $y\in A \setminus P$ alors $x-y=x_0-t+t-y$
(pour n'importe quel $t\in (aB+b) \setminus P$) et comme
$d|(x_0-t)$ et $d|(t-y)$ on a $d|(x_0-y)$. \vskip 2mm \noindent
\tab Ainsi, $d|(x-y)$ pour tout $x,y\in (A \setminus P) \cup
\{x_0\}$ ce qui est en contradiction avec le fait que $(A
\setminus P)\cup \{x_0\}$ soit une base. \vskip 2mm \noindent \tab
Comme $d\hskip -1mm\not \hskip -.13mm |\hskip 1mm a$, il existe un
premier $p$ et un entier $k$ tel que $p^k|d$ et $p^k\hskip
-1mm\not \hskip -.13mm |\hskip 1mm a$. Comme $A\sim
p^kB^{'}+b^{'}$ (puisque presque toute les différences des
éléments de $aB+b$ sont divisibles par $p^k$) et que $A\sim aB+b$
on en déduit que $A\sim paB+c$ pour un certain $c$ ce qui est
impossible puisque $a<pa$ et que $a$ est la raison de $A$. Ceci
prouve donc que $P\subset E$ et par suite que le nombre de parties
essentielles de $A$ est fini, puisque $E$ est fini. \vskip 2mm
\noindent \tab Pour majorer le nombre de parties essentielles en
fonction de $a$, commen\c cons par établir le lemme suivant :
\vskip 2mm \noindent {\bf Lemme 12.---} {\it Soit $A$ une base
additive et $P_1,P_2$ deux essentialités distinctes de $A$ telles
que $P_1\cup P_2\ne A$. Posons
$$d(P_i)=\hbox{\rm pgcd}\{x-y ~|~ x,y\in A \setminus P_i\};\ i=1,2$$
On a $d(P_i)\geq 2$ et $\hbox{\rm pgcd}(d(P_1),d(P_2))=1$.} \vskip
2mm \noindent {\bf Preuve du lemme :} Supposons qu'il existe
$d\geq 2$ tel que $d|d(P_1)$ et $d|d(P_2)$. Considérons un entier
$t\in A \setminus (P_1\cup P_2)$ et la partie $B=A \setminus
(P_1\cap P_2)$. Soit $x\in B$, si $x\notin P_1$ alors
$d|d(P_1)|(x-t)$ et si $x\notin P_2$ alors $d|d(P_2)|(x-t)$. Donc,
pour tout $x,y\in B$, on a $d|(x-y)$ et, par suite, on en déduit
que $B$ ne peut être une base additive. Il s'ensuit, par hypothèse
de minimalité et puisque $P_1\cap P_2\subset P_1,P_2$, que
$P_1=P_1\cap P_2=P_2$ ce qui est, bien sur, absurde par hypothèse.
\fin \vskip 2mm \noindent {\bf Remarque :} L'hypothèse $P_1\cup
P_2\ne A$ est visiblement essentielle pour la preuve de ce lemme,
mais elle est en fait essentielle pour la propriété annoncée. En
effet, si l'on considère $A=\N$, $P_1=2\N$ et $P_2=2\N+1$, on voit
que $P_1$ et $P_2$ sont des essentialités mais que
$d(P_1)=d(P_2)=2$. \vskip 2mm \noindent {\bf Retour à la preuve du
théorème :} Soit $P_1,\cdots,P_s$ les parties essentielles de $A$
et $d_i=d(P_i)\geq 2$ pour tout $i=1,\cdots ,s$. En appliquant le
lemme précédent on voit que pour tout $i=1,\cdots ,s$ on a $A\sim
d_iB_i+b_i$ pour un certain $b_i\in \N$ et un certain $B_i\subset
\N$. Comme les $d_i$ sont premiers entre eux deux à deux, on voit
que l'on a $A\sim d_1\cdots d_sB+b$ pour un certain $b\in \N$ et
un certain $B\subset \N$. Si $a$ désigne la raison de $A$, alors
on a $d_1\cdots d_s|a$ et donc, puisque les $d_i$ sont premiers
entre eux deux à deux, on voit que $a$ est divisible par au moins
$s$ nombres premiers distincts, c'est-à-dire que $s$ est majoré
par la longueur du radical de $a$. \fin \vskip 2mm \noindent {\bf
Remarques :} a) La majoration du nombre de parties essentielles
donnée dans le théorème est en fait la meilleure possible. En
effet, si l'on reprend la suite de bases additives $A_n=p_1\cdots
p_n\N\cup\{p_1\cdots \widehat{p_i}\cdots p_n ~|~ i=1,\cdots ,n\}$
alors le radical de la raison de cette base est $n$ et elle
possède exactement $n$ éléments essentiels qui forment donc
exactement les $n$ parties essentielles de $A_n$. \vskip 2mm
\noindent b) Le fait que les parties essentielles d'une base $A$
soient contenues dans le réservoir $E$ de $A$, permet aussi de
majorer leur nombre par $2^{\# E}$. Toutefois, en termes
d'invariants, il n'est pas intéressant d'utiliser le reservoir
pour majorer le nombre de parties essentielles. Si $A$ et $B$ sont
deux bases telles que $A\sim B$, alors leurs raisons sont égales
et, par suite, le nombre de leurs parties essentielles est
majorable simultanément. Par contre, les réservoirs de $A$ et $B$
peuvent être tr\`es différents et contenir un nombre différent de
parties essentielles. \vskip 2mm \noindent c) Comme pour le cas de
la dessentialistation élémentaire, on peut donner un algorithme de
dessentialisation général : étant donné une base $A$ de parties
essentielles $P_1,\cdots ,P_s$, on note $x_0$ le plus petit
élément de $A \setminus (P_1\cup \cdots \cup P_s)$ et on pose :
$$\begin{array}{lll}
m(A)&=&\displaystyle \hbox{\rm pgcd}\ \left\{x-y ~|~ x,y\in A \setminus (P_1\cup \cdots \cup P_s)\right\}\\
P(A)&=&\displaystyle \left\{\frac{x-x_0}{m(A)} ~|~ x\in A \setminus (P_1\cup \cdots \cup P_s)\right\}\\
\end{array}$$
\tab Le même argument que dans le cas élémentaire montre que $P(A)$ reste une base. On voit que $m(A)=1$ si et seulement si la base $A$ est sans parties essentielles. On considère donc la suite $(P^n(A))_n$ définie par
$$P^0(A)=A,\ \forall n\geq 0\ P^{n+1}(A)=P(P^n(A))$$
\tab Comme pour toute base $B$ on a $m(B)|r(B)$ (o\`u $r(B)$ désigne la raison de $B$) et $r(P(B))=\frac{r(B)}{m(B)}$, on en déduit que pour tout $n\geq 0$ on a $r(P^{n+1}(A))|r(P^n(A))$. Mais l'égalité $r(P^{n+1}(A))=r(P^n(A))$ entraine $m(P^n(A))=1$, ce qui équivaut à $P^n(A)$ sans partie essentielles, ou encore $r(P^n(A))=1$. Ainsi la suite $r(P^n(A))$ est strictement décroissante puis stationnaire égale à $1$ et le rang de stationnarité est majoré par la longueur de l'entier $r(A)$ (i.e la somme des puissances de la décomposition en facteurs premiers de $r(A)$). On en déduit en particulier que si $l$ est la longueur de $r(A)$ alors $P^l(A)$ est égal à la dessentialisée de $A$ modulo une translation et la relation $\sim$.
\vskip 2mm
\noindent
d) On a déjà remarqué qu'il n'existait aucun lien entre la raison et l'ordre d'une base, il faut aussi noter qu'il n'existe aucun lien entre l'ordre et la longueur du radical de la raison d'une base. Ceci implique qu'il est impossible de majorer le nombre de parties essentielles d'une base en fonction de son ordre.
\vskip 2mm
\noindent
\tab Nous finissons cet article en suggérant quelques problèmes ouverts relatifs à l'étude que nous venons de mener.
\vskip 2mm
\noindent
{\bf Problèmes ouverts :}
\vskip 2mm
\noindent
{\it Problème I.---} Nous avons vu que le nombre de parties essentielles n'était pas majorable par une fonction de l'ordre de la base, au contraire du cas des éléments essentiels. Si l'on compte les parties essentielles de cardinal borné, on peut s'interroger sur la possibilité de majorer leur nombre en fonction de $h$. Plus précisément, existe-t-il une fonction $\varphi$ telle que, étant donné $k,h\in \N^*$, pour toute base additive $A$ d'ordre $h$ on ait
$$\#\{P\subset A ~|~ P\ \hbox{partie essentielle de cardinal}\leq k\}\leq \varphi(k,h)\ ?$$
\vskip 1mm \noindent {\it Problème II.---} a) Toute base additive
contient-elle toujours au moins une essentialité? \vskip 2mm
\noindent b) Pour toute base $A$ et toute partie $P\subset A$
telle que $A-P$ ne soit pas une base, existe-t-il $P^{'}\subset P$
telle que $P^{'}$ soit une essentialité de $A$? \vskip 2mm
\noindent c) Est-il vrai qu'une base additive contient toujours
une infinité d'essentialités? \newpage \centerline{\sc
Bibliographie} \vskip 5mm \noindent {\bf [EG]} Paul Erd\"os and
Ronald Graham, {\it On bases with an exact order}, Acta
Arithmetica XXXVII, p. 201-207 (1980). \vskip 2mm \noindent {\bf
[DG]} Bruno Deschamps et Georges Grekos, {\it Majoration du nombre
d'exceptions \`a ce qu'un ensemble de base priv\'e d'un point
reste un ensemble de base}, Journal f\"ur die Reine und Angewandte
Mathematik 539, 45-53 (2001). \vskip 2mm \noindent {\bf [Gr]}
Georges Grekos, {\it Sur l'ordre d'une base additive}, S\'eminaire
de th\'eorie des nombres de Bordeaux 1987-1988, exp. 31. \vskip
2mm \noindent {\bf [MR]} Jean-Pierre Massias et Guy Robin, {\it
Bornes effectives pour certaines fonctions concernant les nombres
premiers}, Journal de théorie des nombres de Bordeaux 8, 215-242
(1996). \vskip 2mm \noindent {\bf [Pl]} Alain Plagne, {\it Sur le
nombre d'éléments exceptionnels d'une base additive}, préprint.

\end{document}